\documentclass[12pt]{article}
\usepackage{amsmath}
\usepackage{amsfonts}
\usepackage{amsthm}
\usepackage{amssymb}
\usepackage{color}
\usepackage{graphicx}
\usepackage{indentfirst}
\usepackage{listings}
\usepackage{mathrsfs}
\usepackage{CJK}
\usepackage{tikz}
\usetikzlibrary{shapes,arrows}

\def\a{\alpha}
\def\b{\beta}
\def\g{\gamma}

\def\e{\varepsilon}

\def\l{\lambda}

\def\n{\nabla}

\def\t{\times}

\def\th{\theta}

\def\i{\infty}

\def\cH{{\cal H}}

\def\cK{{\cal K}}

\def\cP{{\cal P}}

\def\cl{{\cal l}}

\def\dbN{{\mathop{\rm l\negthinspace N}}}

\def\dbR{{\mathop{\rm l\negthinspace R}}}

\def\G{\Gamma}

\def\O{\Omega}
\def\oO{{\overline \O}}

\def\ds{\displaystyle}
\def\ns{\noalign{\smallskip} }

\def\q{\quad}
\def\qq{\qquad}
\def\({\Big (}
\def\){\Big )}
\def\[{\Big[}
\def\]{\Big]}

\newtheorem{lemma}{Lemma}[section]

\newtheorem{remark}{Remark}[section]
\newtheorem{example}{Example}[section]
\newtheorem{theorem}{Theorem}[section]
\newtheorem{corollary}{Corollary}[section]

\newtheorem{definition}{Definition}[section]
\newtheorem{proposition}{Proposition}[section]

\newtheorem{graph}{Graph}[section]

\def\be{\begin{equation*}}
\def\bel{\begin{equation}\label}
\def\ee{\end{equation}}
\def\eee{\end{equation*}}
\def\ba{\begin{array}}
\def\ea{\end{array}}
\def\bt{\begin{theorem}}
\def\et{\end{theorem}}
\def\bc{\begin{corollary}}
\def\ec{\end{corollary}}
\def\br{\begin{remark}}
\def\er{\end{remark}}
\def\bl{\begin{lemma}}
\def\el{\end{lemma}}
\def\bp{\begin{proposition}}
\def\ep{\end{proposition}}
\def\bd{\begin{definition}}
\def\ed{\end{definition}}
\def\bg{\begin{graph}}
\def\eg{\end{graph}}

\makeatletter

\@addtoreset{equation}{section}
\makeatother

\def\q{\quad}
\def\pa{\partial}

\def\wt{\widetilde}
\def\cl{\overline}
\def\={\buildrel \triangle \over =}

\def\sqr#1#2{{\vcenter{\vbox{\hrule height.#2pt
				\hbox{\vrule width.#2pt height#1pt \kern#1pt \vrule width.#2pt}
				\hrule height.#2pt}}}}
\def\signed #1{{\unskip\nobreak\hfil\penalty50
		\hskip2em\hbox{}\nobreak\hfil#1
		\parfillskip=0pt \finalhyphendemerits=0 \par}}
\def\endpf{\signed {$\sqr69$}}

\begin{document}

\title{\bf Carleman estimates for higher order partial differential operators and applications\thanks{This work is partially supported by the NSF of China
		under grants 11971333, 11931011, and by the Science Development Project of Sichuan University under grant 2020SCUNL201. }}

\author{Xiaoyu Fu\thanks{School
        of Mathematics, Sichuan University, Chengdu
        610064, China. {\small\it
            E-mail:} {\small\tt \textcolor{blue}{xiaoyufu@scu.edu.cn}}} \ and Yuan Gao\thanks{School
        of Mathematics, Sichuan University, Chengdu
        610064, China. {\small\it
            E-mail:} {\small\tt \textcolor{blue}{scumathgaoyuan@163.com}}}}

\date{}

\maketitle

\begin{abstract}

In this paper, we obtain a  Carleman estimate for the higher order partial differential operator $\cP\=\a\pa_t+\pa_x^n$ (with $\a\in\mathbb{R}, n\in\mathbb{Z}_{\ge2})$.
In the process of establishing this estimate, we developed a new method,
which is called the ``back-propagation method'' (the BPM, for short). This method can also be used to
build up Carleman estimates for some other partial differential operators, and might provide assistance
with corresponding numerical analyses. As an application of the above-mentioned Carleman estimate,
we proved the conditional stability of a Cauchy problem for a time fractional diffusion equation with $\frac{1}{3}$-order.

\end{abstract}

\bigskip
\noindent{\bf AMS Subject Classifications:}.  34K37; 35R30; 35R25
\bigskip

\noindent{\bf Key Words}. higher order partial differential operator, Carleman estimate, back-propagation method.

\section{Introduction and main result }\label{s1}

In 1939,  Carleman \cite{c} showed that a second order elliptic operator in dimensional two enjoys  the unique continuation property. The technique used there is called ``Carleman weight inequality"  and has become one of the major tools in the study of the unique continuation property, control problems and inverse problems for partial differential equations.

The Carleman estimate can be regarded as a weighted energy inequality, which can be described as follows:

Let $\Omega$ be a
connected open set in $\mathbb{R}^{n}$, and let $P=P(x,D)$ be a differential operator of order $n$ in $\Omega$. Assume that there is a suitable function
$\phi(\cdot)\in C^\i(\oO;\mathbb{R})$ satisfying $\ds \n\phi(x)\neq 0$, $x\in\O$. Let $\th=e^{\l\phi}$. We say that the Carleman type estimate holds for $P$ if there exists a constant $C>0$ such that
$$\ba{ll}\ds
\sum_{0\le m<n}\l^{2(n-m)-1}\int_{\Omega}\th^2|D^{m}v|^{2}dx\leq C\int_{\Omega}|\th P(x,D)v|^{2}dx,
\ea
$$
where $v\in C_0^\i(\O)$ and  $\l>0$ is a parameter.

Up to now, there are numerous results on Carleman estimates for the second order partial differential operators,
the corresponding applications are well understood (for example, see \cite{FLZ} and the references cited therein).
For Carleman estimates of  higher order partial differential operators,  we refer to  \cite{CFT, G, GK, XCY, Zheng, ZZ}  for the fourth order  parabolic-type operators and fourth order Schr\"odinger operators, \cite{JF}  for the sixth-order parabolic operators. Based on an identity of Treves (see  \cite[Lemma 17.2.2]{Hormander}), \cite{HTX, HHZ} obtained unique continuation properties for higher Order parabolic equations and Schr\"odinger Equations. Recently \cite{Isaza} established a Carleman estimate for  high order equations of Korteweg-de Vries type with the weight function $\th=e^{\l x}$ in one dimensional case. Compared to Carleman estimates for second order partial differential operators, the computation of Carleman estimates for higher order partial differential operators is much more complicated.

In this paper, we aim at establishing a Carleman estimate for the operator
$$\cP =\a\pa_t+ \pa_x^n,\ \a\in\mathbb{R}\backslash\{0\},\ n\in\mathbb{Z}_{\ge2},$$ with an exponential-type weight function $\th(t, x)$ (which  will be given later).  More importantly, using $\cP$ as a carrier, we will introduce the ``backpropagation method'' (the BPM, for short) in the building of the Carleman estimate for $\cP$. \par
Next, we will explain what role the BPM will play in the process of establishing the Carleman
estimate for $\cP$. As we know, elementary calculus can be enough to grasp the main idea of Carleman estimate (See \cite[Chapter 1]{FLZ}). To obtain the core, there is an important ingredient which should be emphasized, i.e.,
     the principal part of our operator $\th (\cP \th^{-1})$. Noting that the principal part of our operator contains the $n$-th derivative with respect to the $x$ variable,
hence the decomposition and the computation will cause many difficulties.  For the reader's convenience, we simply choose  the following  weight functions:

\bel{psi-0}
 \psi(t,x)=(x-x_{0})^2-\b(t-t_0)^2,\q  \ell(t,x)=\l\psi(t,x),\q \th(t,x)=e^{\ell(t,x)}
  \ee
 where $x_0,t_0\in\mathbb{R}$, $\b\in\mathbb{R}_{+}$.

For $v(t,x)\in C^{\i}(\mathbb{R}\times\mathbb{R};\mathbb{R})$, set
 $
w=\th v=e^{\ell} v
$ with $\ell$ given by (\ref{psi-0}). Then
\bel{as00}\ba{ll}\ds
\th\cP v&\ds=\a\th \pa_t v+\th\pa^n_x v\\
\ns&\ds=\a w_t-\a\ell_t w+I_1(w)+I_2(w)+I_3(w).
\ea\ee
Here $I_j(w), j=1,2,3$ are chosen in the following manner: $$\th\pa_x^n v=I_1(w)+I_2(w)+I_3(w).$$

 Note that the order of $\l$ equals to the sum of orders of $\pa_x^j\ell$, $j\ge1$. Then we decompose $\th\pa_x^nv$ as following:

\begin{itemize}

\item The item ${I_1(w)}$ is the sum of all such terms that are the highest order terms of $\l$ with an odd
order derivative of $x$, and the second highest orders terms of $\l$ with an even order derivative of $x$;

 \item The item ${I_2(w)}$ is the sum of all such terms that are the highest order terms of $\l$ with an even
order derivative of $x$, and the second highest orders terms of $\l$ with an odd order derivative of $x$;

\item The item ${I_3(w)}$ is the sum of all other terms in $\th\pa_x^n v$.
  \end{itemize}
It deserves mentioning that ${I_3(w)}$ consists of only lower order terms, compared with ${I_1(w)}$ and ${I_2(w)}$.

 By (\ref{as00}), we have

\bel{1128-a0}
\th\cP vI_2=\a w_t I_2+I_1I_2+|I_2|^2+I_2(I_3-\a\ell_t w).
 \ee

From (\ref{1128-a0}), it is easy to see that:
\bel{23061404}\ba{ll}\ds
\th^2|\cP v|^2\ge 2I_1I_2+2\a w_tI_2-|{I_3(w)}-\a\ell_t w|^2.
\ea\ee
Comparing ${I_3(w)}-\a\ell_t w$ with $I_1, I_2$, we know $|{I_3(w)}-\a\ell_t w|^2$ contains only lower order term.

Our keys to establishing the Carleman estimate for $\cP$ are as follows:

\begin{itemize}
\item The decomposition of principal operator $\th \pa^n_x v$ (see Proposition \ref{lmSW01}).
\item The choice of $I_j(w) (j=1, 2, 3)$ (see Proposition \ref{lmSW02}).
\item The estimates of $I_1(w)I_2(w)$ and $w_tI_2(w)$(see Propositions \ref{22051602} and \ref{23052501}).
\end{itemize}

The BPM developed in this paper is used to solve the third key above. More precisely, it can help
us to simply prove that the highest order energy terms in $I_1(w)I_2(w)$ have the positive sign, while the
``bad'' terms yielded from $w_tI_2(w)$, such as $w_tw_x$ and $w_{tx}w_{xx}$, are indeed lower order terms.
Our BPM is inspired by ``back-propagation'', which is widely used in the field of machine learning, dating back to \cite{DER} published in 1986. In machine learning, it is difficult to calculate derivatives forward, so people adopt a backward method, which was proposed in the 1970s, to solve differential problems of nested functions. Now people working on machine learning will draw a network graph to calculate a derivative backward, which can be viewed as the embryo of our BPM.

Throughout of this paper, we use $C_n^k$ to denote combinatorial numbers, $v_t$ or $\pa _t v$ to represent the derivative of $v$ in the time variable, and $v_x$ or $\pa _x v$ to denote the derivative of $v$ in the space variable. In what follows,  we will use $C$ to denote a generic positive constant which may vary from line to line. For $k\in\mathbb{Z}_{\ge1}$, we denote by $O(\l^k)$ a function of order $\l^k$ for large $\l$.  We use ``$(\cdot)$'' or ``etc.'' to denote  such terms whose concrete forms does not have to be given.

 The main result of this paper is as follows:

\bt\label{thmSW00}
Let  $T>0$,\ $L>0$,  $\a\in\mathbb{R}\backslash\{0\},\ n\in\mathbb{Z}_{\ge2}$.  Let $\psi=(x-x_0)^2-\beta(t-t_0)^2$ with $x_0>L$, $0<t_0<T$, $\beta\in\mathbb{R}_+$. For any $v \in C_0^\i([0,T]\t[0,L];\mathbb{R})$, and $w=e^{\l\psi}v$, the following inequality holds:
\bel{22020801}\ba{ll}\ds
\sum^{n-1}_{m=0}\int_{0}^T\int_0^L\[n^2C_{n-1}^{m}\l^{2n-2m-1}\psi_x^{2n-2m-2}\psi_{xx}+O(\l^{2n-2m-3})\]|\pa_x^{m}w|^2dxdt\\
\ns\ds\le  \int_{0}^T\int_0^Le^{2\l\psi}\big|\a \pa_t v+\th\pa^n_x v\big|^2dxdt.\ea\ee
\et

\br
The main innovation of this paper is the BPM. It is not only applicable to different
partial differential operators, but also to different forms of the weight function $\varphi$. In addition, it can be
used for the estimation of boundary terms. These will be given in our forthcoming work.
\er

\br
The inequality (\ref{22020801}), with $\a=-1$ and $n=4$, was built up in \cite{XCY} by a different
way from ours. Based on it, one can obtain the conditional stability in a Cauchy problem for a half-order
fractional diffusion equation.
\er

\br
The reason why we do not convert the function $w$ on the left side of (\ref{22020801})  into $v$ is as follows: This way not only gives us the explicit coefficients of all energy terms, but also simplifies
the computation. In fact, this transformation can be easily implemented in the following way: Since $w=e^{\l\psi}  v$, we can find a constant $C_1>0$ such that
 \bel{1107-c}\ba{ll}\ds
 \sum^{n-1}_{m=0}\int_{0}^T\int_0^L\l^{2n-2m-1}e^{2\l\psi}|\pa_x^{m}v|^2dxdt\\
 \ns\ds
 =\sum^{n-1}_{m=0}\int_{0}^T\int_0^L\l^{2n-2m-1}e^{2\l\psi}|\pa_x^{m}(e^{-\l\psi}w)|^2dxdt\\
 \ns\ds\le C_1\sum^{n-1}_{m=0}\int_{0}^T\int_0^L\l^{2n-2m-1}|\pa_x^{m}w|^2dxdt.
 \ea\ee
Then, the combination of  (\ref{22020801}) and (\ref{1107-c}) yields  that there exists a $\l_0>0$, such that for any $\l>\l_0$, we have
 \bel{22020-09}\ba{ll}\ds
\sum^{n-1}_{m=0}\int_{0}^T\int_0^L\l^{2n-2m-1}e^{2\l\psi}|\pa_x^{m}v|^2dxdt\le C\int_{0}^T\int_0^Le^{2\l\psi}|\cP v|^2dxdt\ea\ee
for some $C>0$.
 \er

The rest of this paper is organized as follows. In section 2, we give some preliminaries. The estimations of $I_1I_2$ and $w_tI_2$ are given in sections 3 and 4, respectively. In section 5, we give our proof of Theorem \ref{thmSW00}. As its application, we consider an inverse problem for a  time-fractional diffusion equation in section 6.

\section{ Some Preliminaries }

This section presents preliminaries.  We start with introducing notation.

For $m_1, m_2\in \mathbb{R}$, denote
\be
\sum_{j={m_1}}^{m_2} a_j\=\sum_{j\in [m_1,m_2]\cap\mathbb{Z}} a_j,\qq \prod_{j={m_1}}^{m_2} a_j\=\prod_{j\in [m_1,m_{m_2}]\cap\mathbb{Z}} a_j,
\eee
with the convention that $\ds \sum_{j\in\varnothing}a_j\=0$ and $\ds \prod_{j\in\varnothing}a_j\=1$.\par

For $j,k\in \mathbb{Z}$, denote
\be
C_{j}^{k}\=\left\{\ba{ll}
\ns\ds \frac{j!}{k!(j-k)!}, & \mbox{ if } 0\le k\le j, \\
\ns\ds 0, & \mbox{ else }.
\ea\right.\eee

\subsection{ An Identity about Combinatorial Numbers}

\bp\label{22051601}
For $(n,m)\in\{\mathbb{Z}^2|n\ge2,\ 0\le m\le n-1\}$, define
\bel{1008-0}\ba{ll}\ds
\cK_{n,m}&\ds\=-C_{n}^{2}C_{n}^{m}C_{n-2}^{m}-C_{n}^2\sum_{k=1}^m(-1)^k(C_{n}^{m+k}C_{n-2}^{m-k}+C_{n}^{m-k}C_{n-2}^{m+k}) \\
\ns&\ds\q+\frac{2n-2m-1}{2}\sum_{k=0}^m(-1)^k(1+2k)(C_{n}^{m+1+k}C_{n}^{m-k}).
\ea\ee
Then
\bel{1008-00}
\cK_{n,m}=\frac{n^2}{2}C_{n-1}^{m}.
\ee
\ep

\br
In (\ref{1008-0}), the form of  $\cK_{n,m}$ is the highest order coefficients of energy terms involved in the estimation of $I_1I_2$. We set this Proposition because it is not an obvious result.
\er

\ \\
{\bf Proof of Proposition \ref{22051601}. }  For $(n,m)\in\{\mathbb{Z}^2|n\ge2,\ 0\le m\le n-1\}$, we simply write
 \bel{1008-1}
 \cK_{n,m}=C_{n}^{2}\cH_1+\frac{2n-2m-1}{2}\cH_2,
 \ee
where
 \bel{1008-2}\left\{\ba{ll}\ds
\cH_1= -C_{n}^{m}C_{n-2}^{m}-\sum_{k=1}^m(-1)^k(C_{n}^{m+k}C_{n-2}^{m-k}+C_{n}^{m-k}C_{n-2}^{m+k}),\\
\ns\ds\cH_2=\sum_{k=0}^m(-1)^k(1+2k)(C_{n}^{m+1+k}C_{n}^{m-k}).
 \ea\right.\ee

We first calculate $\cH_1$. By replacing the index of summation, one can get:
\bel{22051603}\ba{ll}
\cH_1&\ds= -C_{n}^{m}C_{n-2}^{m}-\sum_{k=1}^m(-1)^k(C_{n}^{m+k}C_{n-2}^{m-k}+C_{n}^{m-k}C_{n-2}^{m+k})\\
\ns&\ds =(-1)^{m+1}\sum_{q=m}^m(-1)^{q}C_{n}^{q}C_{n-2}^{2m-q}+(-1)^{m+1}\sum_{q=m+1}^{2m}(-1)^{q}C_{n}^{q}C_{n-2}^{2m-q} \\
\ns&\ds\q +(-1)^{m+1}\sum_{q=0}^{m-1}(-1)^{q}C_{n}^{q}C_{n-2}^{2m-q} \\
\ns&\ds =(-1)^{m+1}\sum_{q=0}^{2m}(-1)^{q}C_{n}^{q}C_{n-2}^{2m-q}.
\ea\ee
Meanwhile, we have the identity:
\bel{22051604}
(1-x)^n(1+x)^{n-2}=(1-x^2)^{n-2}(x^2-2x+1), \q x\in \dbR,\  n\in\dbN.
\ee
Comparing the coefficients of $x^{2m}$ (with $m\in\dbN$) on the both sides of (\ref{22051604}), we find
\bel{22051605}
\ds \sum_{q=0}^{2m}(-1)^{q}C_{n}^{q}C_{n-2}^{2m-q}=C_{n-2}^{m}(-1)^m+C_{n-2}^{m-1}(-1)^{m-1}, \q m\in\dbN.
\ee
Combining (\ref{22051603}) and (\ref{22051605}), we have
 \bel{1008-0a}
 \cH_1=C_{n-2}^{m-1}-C_{n-2}^{m}.
 \ee

Next, let us calculate $\cH_2$. By a similar way dealing with $\cH_1$, we can obtain
\bel{22051606}\ba{ll}\ds
\cH_2&\ds= \sum_{k=0}^m(-1)^k(1+2k)(C_{n}^{m+1+k}C_{n}^{m-k})\\
\ns&\ds =\sum_{k=0}^m\[(-1)^k(m+1+k)+(-1)^{k+1}(m-k)\]C_{n}^{m+1+k}C_{n}^{m-k}\\
\ns&\ds =\sum_{q=m+1}^{2m+1}(-1)^{q-m-1}qC_{n}^{q}C_{n}^{2m+1-q}+\sum_{q=0}^m(-1)^{m-q-1}qC_{n}^{2m+1-q}C_{n}^{q} \\
\ns&\ds =(-1)^{m+1}\sum_{q=1}^{2m+1}(-1)^{q}qC_{n}^{q}C_{n}^{2m+1-q}.
\ea\ee

At the same time, we have the identity:
\bel{22051607}
(1+x)^n\frac{d}{d x}(1-x)^n=-n(1+x)(1-x^2)^{n-1}, \q x\in \dbR,\  n\in\dbN.
\ee
Comparing  the coefficient of $x^{2m}$ (with $m\in\dbN$) on the both sides of (\ref{22051607}), we find:
\bel{22051608}\ba{ll}\ds
\q \sum_{k=0}^{2m}C_{n}^{2m-k}(-1)^{k+1}(k+1)C_{n}^{k+1}=\sum_{q=1}^{2m+1}(-1)^{q}qC_{n}^{q}C_{n}^{2m+1-q}=-n(-1)^{m}C_{n-1}^m.
\ea\ee
Combining (\ref{22051606}) and (\ref{22051608}), we have
 \bel{1008-0b}
 \cH_2=nC_{n-1}^m.
 \ee

Finally, it follows from (\ref{1008-1}), (\ref{1008-0a}) and (\ref{1008-0b}) that
\bel{1008-0c}\ba{ll}
\cK_{n,m}&\ds= C_{n}^{2}(C_{n-2}^{m-1}-C_{n-2}^{m})+\frac{n(2n-2m-1)}{2}C_{n-1}^{m} \\
\ns&\ds =\frac{n^2}{2}C_{n-1}^{m}.
\ea\ee
This completes the proof of Proposition \ref{22051601}. \endpf

\subsection{ Decomposition of $\th \pa_x^n v$}

In this subsection, we introduce a decomposition of $\th \pa_x^n v$, which will play an important role in the proof of our main result. \par

\bp\label{lmSW01} Let $\th$ be given by (\ref{psi-0}). Set $w=\th v$. Then

\bel{22030802}
\ds\th\pa_x^{n} v=\sum_{r+2s+m=n}\[(-1)^{s}\frac{1}{s!}(\prod_{l=0}^{s-1}{C_{n-2l}^2})\ell_{xx}^s\]\[(-1)^r C_{r+m}^m\ell_x^r\pa_x^m v\],
\ee
where $r,s,m\in \mathbb{Z}_{\ge0}$.
\ep

\br
For any fixed $n$, the decomposition of $\th\pa_x^n w$ can be obtained by virtue of an iteration. Indeed, it is clear that
$$\th\pa_x^0 v=w,\q \th\pa_x^1 v=\th\pa_x^1(\th^{-1}w)=\pa_x w-\ell_x w.$$
Assuming $f=\th\pa_x^n v$, noting that $\th_x=\th \ell_x$, we have
\bel{22020802}\ba{ll}
\ns\th\pa_x^{n+1} v&\ds=\th\pa_x\(\frac{f}{\th}\)=\th\(\frac{f_x\th-\th_x f}{\th^2}\)=\pa_x f-\ell_x f\\
\ns&\ds=\pa_x(\th\pa_x^n v)-\ell_x(\th\pa_x^n v).
\ea\ee
\er

As we explained before, the  order of $\l$ equals to the sum of orders of $\pa_x^j\ell$, $j\ge1$. In fact, we only care about terms with $s=0,1$ in $\th\pa_x^{n} v$. The terms involved in (\ref{22030802}) with $s\ge 2$ only yield lower order terms. Truncating $\theta\pa_x^nv$ at $s=0,1$, we have
\bel{23091501}\ba{ll}
\theta\pa_x^{n}v
&\ds=\sum_{r+m=n}(-1)^{r}C_{n}^m\ell_x^r\pa_x^mw-\sum_{r+m=n-2}(-1)^{r}C_{n}^2C_{n-2}^m\ell_x^r\ell_{xx}\pa_x^mw+\cdots. \\
\ns &\ds=\sum_{m=0}^{n}(-1)^{n-m}C_n^m\ell_x^{n-m}\pa_x^mw\\
\ns&\ds\q-C_n^2\ell_{xx}\sum_{m=0}^{n-2}(-1)^{n-2-m}C_{n-2}^m\ell_x^{n-2-m}\pa_x^mw
+\cdots.
\ea\ee

\bp\label{lmSW02}
Under the assumption of Proposition \ref{lmSW01}, we have
\bel{0425-a0}
\th \pa_x^n v=I_1(w)+I_2(w)+I_3(w),
\ee
where
\bel{1128-b0}\left\{\ba{ll}
\ns\ds I_1(w)=\sum_{odd\ j\in[0,n]}C_{n}^{j}(-1)^{n-j}\ell_x^{n-j}\pa_x^{j}w \\
\ns\ds\q\q\qq -C_{n}^2\ell_{xx}\sum_{even\ k\in[0,n-2]}C_{n-2}^{k}(-1)^{n-2-k}\ell_x^{n-2-k}\pa_x^{k}w, \\
\ns\ds I_2(w)=\sum_{even\ k\in[0,n]}C_{n}^{k}(-1)^{n-k}\ell_x^{n-k}\pa_x^{k}w \\
\ns\ds\q\q\qq -C_{n}^2\ell_{xx}\sum_{odd\ j\in[0,n-2]}C_{n-2}^{j}(-1)^{n-2-j}\ell_x^{n-2-j}\pa_x^{j}w,\\
\ns\ds I_3(w)\=\th \pa_x^n v-I_1(w)-I_2(w).
\ea\right.\ee
\ep

For concrete $n$, by Proposition \ref{lmSW02}, we can obtain the decomposition of $\th\pa_x^n w$ immediately.

\begin{example}\label{I1I2n4} In the case that $n=4$, we have
\bel{22053001}\left\{\ba{ll}\ds
I_1=
-4\ell_x\pa_x^{3}w
-4\ell_x^{3}\pa_xw
-6\ell_{xx}\pa_x^{2}w
-6\ell_x^{2}\ell_{xx}w, \\
\ns\ds
I_2=
\pa_x^{4}w
+6\ell_x^{2}\pa_x^{2}w
+\ell_x^{4}w
+12\ell_x\ell_{xx}\pa_xw. \\
\ea\right.\ee
\end{example}

\begin{example}\label{I1I2n5} In the case that $n=5$, we have
\be\left\{\ba{ll}\ds
I_1=
\pa_x^{5}w
+10\ell_{x}^{2}\pa_x^{3}w
+5\ell_{x}^{4}\pa_xw
+30\ell_{x}\ell_{xx}\pa_x^{2}w
+10\ell_{x}^{3}\ell_{xx}w, \\
\ns\ds
I_2=
-5\ell_{x}\pa_x^{4}w
-10\ell_{x}^{3}\pa_x^{2}w
-\ell_{x}^{5}w
-10\ell_{xx}\pa_x^{3}w
-30\ell_{x}^{2}\ell_{xx}\pa_xw. \\
\ea\right.\eee
\end{example}

Our proof of Proposition \ref{lmSW01} is based on the iteration (\ref{22020802}). The following example will help us to understand our idea.

\begin{example}\label{23091502} If we have already known that $\theta\pa_x^4v=12\ell_x\ell_{xx}\pa_xw-4\ell_x^3\pa_xw-6\ell_x^2\ell_{xx}w+\cdots$ (see Example \ref{I1I2n4}), we can find all the coefficients in $\ell_x^r\ell_{xx}^s\pa_x^mw$ with $r+2s+m=5$ in $\theta\pa_x^5v$. Let us see how to obtain the coefficient of $\ell_x^2\ell_{xx}\pa_xw$ in $\theta\pa_x^5v$.

By (\ref{22020802}), one knows that $\theta\pa_x^5v=\pa_x(\th\pa_x^4 v)-\ell_x(\th\pa_x^4 v)$. So one only needs to find all terms $A$ in $\th\pa_x^4 v$ such that
$$\pa_xA=c\cdot\ell_x^2\ell_{xx}\pa_xw+\cdots, \mbox{ or } -\ell_x A=c\cdot\ell_x^2\ell_{xx}\pa_xw+\cdots,$$
where $c\in\mathbb{R}$ is the coefficient.
There are $3$ terms one can find, which are
\bel{1128-0}\left\{\ba{ll}
\ns\ds -\ell_x\cdot(12\ell_x\ell_{xx}\pa_xw)=-12\ell_x^2\ell_{xx}\pa_xw, \\
\ns\ds \pa_x(-4\ell_x^3\pa_xw)=-12\ell_x^2\ell_{xx}\pa_xw+\cdots, \\
\ns\ds \pa_x(-6\ell_x^2\ell_{xx}w)=-6\ell_x^2\ell_{xx}\pa_xw+\cdots.
\ea\right.\ee
Adding them up leads to the coefficient $-30$. This  matches what we got in Example \ref{I1I2n5}.
\end{example}

\ \\
{\bf Proof of Proposition \ref{lmSW01}.}
It is easy to check that
\be\ba{ll}\ds
\th\pa_x v=w_x-\ell_xw, \\
\ns\ds\th\pa_x^2 v=w_{xx}-2\ell_x w_x+\ell_x^2 w-\ell_{xx}w
\ea\eee \\
satisfy (\ref{22030802}). We inductively assume that (\ref{22030802}) holds for some $n\ge 2$. We will prove that (\ref{22030802}) holds  for $n+1$.

By (\ref{22020802}), we know that  $\th\pa_x^{n+1} v=\pa_x(\th\pa_x^{n} v)-\ell_x(\th\pa_x^{n} v)$. Considering a general term $\ell_x^r\ell_{xx}^s\pa_x^mw$ in $\th\pa_x^{n+1}v$ with $r+2s+m=n+1$, proceeding exactly the same analysis of (\ref{1128-0}), we know that there are $3$ terms in $\th\pa_x^{n} v$ contributing to $\ell_x^r\ell_{xx}^s\pa_x^m$:
\bel{23061401}\left\{\ba{ll}
\ns\ds (-\ell_x)\cdot c_1\ell_x^{r-1}\ell_{xx}^s\pa_x^mw=-c_1 \ell_x^r\ell_{xx}^s\pa_x^mw, \\
\ns\ds \pa_x(c_2\ell_x^{r+1}\ell_{xx}^{s-1}\pa_x^mw)=(r+1)c_2 \ell_x^r\ell_{xx}^s\pa_x^mw+\cdots, \\
\ns\ds \pa_x(c_3\ell_x^r\ell_{xx}^s\pa_x^{m-1}w)=c_3 \ell_x^r\ell_{xx}^s\pa_x^mw+\cdots,
\ea\right.\ee
where
\be\left\{\ba{ll}
\ns\ds c_1=(-1)^{r-1+s}\frac{1}{s!}\(\prod_{l=0}^{s-1}{C_{n-2l}^2}\)C_{r-1+m}^m,\\
\ns\ds c_2=(-1)^{r+1+s-1}\frac{1}{(s-1)!}\(\prod_{l=0}^{(s-1)-1}{C_{n-2l}^2}\)C_{r+1+m}^m, \\
\ns\ds c_3=(-1)^{r+s}\frac{1}{s!}\(\prod_{l=0}^{s-1}{C_{n-2l}^2}\)C_{r+m-1}^{m-1}.
\ea\right.\eee
Adding 3 formulas in (\ref{23061401}) up, noting that $r+2s+m=n+1$, we can get the coefficient of $\ell_x^r\ell_{xx}^s\pa_x^mw$ in $\th\pa_x^{n+1} v$. Considering  whether $r,s,m$ equal(s) zero, we discuss the following 7 cases:

1. $r\neq0,\ s\neq0, m\neq0:$
\bel{23110302}\ba{ll}
\ns\ds d&\ds=-c_1+(r+1)c_2+c_3\\
\ns&\ds=(-1)^{r+s}\frac{1}{s!}\(\prod_{l=0}^{s-1}{C_{n-2l}^2}\)\[C_{r-1+m}^m+\frac{(r+1)s}{C_{n-2(s-1)}^2}C_{r+1+m}^m+C_{r+m-1}^{m-1}\]\\
\ns&\ds =(-1)^{r+s}\frac{1}{s!}(\prod_{l=0}^{s-1}C_{n-2l}^2)\[C_{r-1+m}^m+\frac{(r+1)s}{C_{r+m+1}^2}C_{r+1+m}^m+C_{r+m-1}^{m-1}\] \\
\ns&\ds=(-1)^{r+s}\frac{1}{s!}(\prod_{l=0}^{s-1}C_{n-2l}^2)\frac{r+2s+m}{r+m}C_{r+m}^{m} \\
\ns&\ds =(-1)^{r+s}\frac{1}{s!}(\prod_{l=0}^{s-1}C_{n-2l}^2)\frac{n+1}{n+1-2s}C_{r+m}^{m} \\
\ns&\ds =(-1)^{r+s}\frac{1}{s!}(\prod_{l=0}^{s-1}C_{n+1-2l}^2)C_{r+m}^{m}.
\ea\ee

2. $r=0,\ s\neq0, m\neq0:\ c_1=0$,
\be\ba{ll}
\ns\ds d&\ds=(r+1)c_2+c_3\\
\ns&\ds=(-1)^{s}\frac{1}{s!}\(\prod_{l=0}^{s-1}{C_{n-2l}^2}\)\[\frac{s}{C_{n-2(s-1)}^2}C_{1+m}^m+C_{m-1}^{m-1}\]\\
\ns&\ds =(-1)^{s}\frac{1}{s!}(\prod_{l=0}^{s-1}C_{n+1-2l}^2)\frac{n+1-2s}{n+1}\[\frac{2s}{m}+1\] \\
\ns&\ds=(-1)^{s}\frac{1}{s!}(\prod_{l=0}^{s-1}C_{n+1-2l}^2)\frac{n+1-2s}{n+1}\ \frac{2s+m}{m} \\
\ns&\ds =(-1)^{s}\frac{1}{s!}(\prod_{l=0}^{s-1}C_{n+1-2l}^2).
\ea\eee

3. $r\neq0,\ s=0, m\neq0:\ c_2=0$,
\be\ba{ll}
\ns\ds d&\ds=-c_1+c_3=C_{r-1+m}^m+C_{r+m-1}^{m-1}=C_{r+m}^m.
\ea\eee

4. $r\neq0,\ s\neq0, m=0:\ c_3=0$,
\be\ba{ll}
\ns\ds d&\ds=-c_1+(r+1)c_2\\
\ns&\ds =(-1)^{s}\frac{1}{s!}\(\prod_{l=0}^{s-1}{C_{n-2l}^2}\)\[C_{r-1+0}^0+\frac{(r+1)s}{C_{n-2(s-1)}^2}C_{r+1+0}^0\]\\
\ns&\ds =(-1)^{s}\frac{1}{s!}\(\prod_{l=0}^{s-1}{C_{n-2l}^2}\)\[1+\frac{(r+1)s}{C_{n-2(s-1)}^2}\]\\
\ns&\ds =(-1)^{s}\frac{1}{s!}\(\prod_{l=0}^{s-1}C_{n+1-2l}^2\)\frac{n+1-2s}{n+1}\[1+\frac{2(r+1)s}{(n-2s+2)(n-2s+1)}\]\\
\ns&\ds =(-1)^{s}\frac{1}{s!}\(\prod_{l=0}^{s-1}C_{n+1-2l}^2\)\frac{r}{n+1}\[1+\frac{(r+1)(n+1-r)}{(r+1)r}\]\\
\ns&\ds =(-1)^{s}\frac{1}{s!}\(\prod_{l=0}^{s-1}C_{n+1-2l}^2\).
\ea\eee

5, $r=0,\ s=0,\ m\neq0:\ c_1=c_2=0$,
\be\ba{ll}
\ns\ds d&\ds=c_3=C_{m-1}^{m-1}=1.\\
\ea\eee

6. $r=0,\ s\neq0,\ m=0:\ c_1=c_3=0$,
\be\ba{ll}
\ns\ds d&\ds=(r+1)c_2 \\
\ns&\ds =(-1)^{s}\frac{1}{(s-1)!}\(\prod_{l=0}^{s-2}{C_{n-2l}^2}\)C_{1}^0 \\
\ns&\ds =(-1)^{s}\frac{1}{s!}\frac{n+1}{2}\frac{n(n-1)}{2}\frac{(n-2)(n-3)}{2}...\frac{3\times2}{2} \\
\ns&\ds =(-1)^{s}\frac{1}{s!}\frac{(n+1)n}{2}...\frac{4\times3}{2} \\
\ns&\ds =(-1)^{s}\frac{1}{s!}(\prod_{l=0}^{s-1}C_{n+1-2l}^2).
\ea\eee

7. $r\neq0,\ m=0,\ s=0:\ c_2=c_3=0$,
\be\ba{ll}
\ns\ds d&\ds=-c_1=(-1)^{r}\frac{1}{s!}\(\prod_{l=0}^{s-1}{C_{n-2l}^2}\)C_{r-1}^0=(-1)^{r}\frac{1}{s!}\(\prod_{l=0}^{s-1}{C_{n-2l}^2}\).
\ea\eee

We note that the results of cases 2-7 are consistent with those of case 1. Then we complete the proof of Proposition \ref{lmSW01}.
\endpf

\subsection{Introduction of Back-Propagation Method}

In this subsection, we will show how the back-propagation method works in this paper.

Let $(A,B)\in\mathbb{Z}_{\ge0}\times\mathbb{Z}_{\ge0}$ be a pair of positional parameters (which can be replaced), and $F(A,B)$ be a given term located at position $(A, B)$. For $(X, Y)\in\mathbb{Z}\t\mathbb{Z}$, we define
\bel{24011801}
F(X,Y)=\left\{\begin{matrix}
\ns 0,\ & X\notin[0,A]\mbox{ or }Y\notin[0,B], \\
\ns r_{XY}F(X-1,Y)+s_{XY}F(X,Y-1),\ & else, \\
\end{matrix}\right.\ee
where $r_{XY},s_{XY}\in\mathbb{R}$ are known.

The relationship given in (\ref{24011801}) can be represented by the following propagation graph:
\bg\label{graph1}
\begin{CJK*}{GBK}{amsmath}
\tikzstyle{process} = [rectangle,rounded corners, minimum width=1.2cm,minimum height=0.7cm,text centered,text width =1.4cm, draw=black,fill=white]
\tikzstyle{arrow1} = [thick,->,>=stealth]
{$\\ \mbox{ }\q\q\q$}
\begin{tikzpicture}
\node (process1) [process]{$F(0,0)$};
\node (process2) [process,right of=process1,xshift=1.8cm] {$F(0,1)$};
\node (process3) [process,right of=process2,xshift=1.8cm] {...};
\node (process4) [process,right of=process3,xshift=1.8cm] {$F(0,B)$};
\node (process5) [process,below of=process1,yshift=-.2cm]{$F(1,0)$};
\node (process6) [process,right of=process5,xshift=1.8cm] {$F(1,1)$};
\node (process7) [process,right of=process6,xshift=1.8cm] {...};
\node (process8) [process,right of=process7,xshift=1.8cm] {$F(1,B)$};
\node (process9) [process,below of=process5,yshift=-.2cm]{...};
\node (process10) [process,right of=process9,xshift=1.8cm] {...};
\node (process11) [process,right of=process10,xshift=1.8cm] {...};
\node (process12) [process,below of=process8,yshift=-.2cm] {...};
\node (process13) [process,below of=process9,yshift=-.2cm]{$F(A,0)$};
\node (process14) [process,right of=process13,xshift=1.8cm] {$F(A,1)$};
\node (process15) [process,right of=process14,xshift=1.8cm] {...};
\node (process16) [process,right of=process15,xshift=1.8cm] {$F(A,B)$};
\draw [arrow1] (process2) -- node[anchor=south]{$s_{01}$} (process1);
\draw [arrow1] (process3) -- node[anchor=south]{$s_{02}$} (process2);
\draw [arrow1] (process4) -- node[anchor=south]{$s_{0B}$} (process3);
\draw [arrow1] (process6) -- node[anchor=south]{$s_{11}$} (process5);
\draw [arrow1] (process7) -- node[anchor=south]{$s_{12}$} (process6);
\draw [arrow1] (process8) -- node[anchor=south]{$s_{1B}$} (process7);
\draw [arrow1] (process10) -- node[anchor=south]{$s_{\cdot 1}$} (process9);
\draw [arrow1] (process11) -- node[anchor=south]{$s_{\cdot 2}$} (process10);
\draw [arrow1] (process12) -- node[anchor=south]{$s_{\cdot B}$} (process11);
\draw [arrow1] (process14) -- node[anchor=south]{$s_{A1}$} (process13);
\draw [arrow1] (process15) -- node[anchor=south]{$s_{A2}$} (process14);
\draw [arrow1] (process16) -- node[anchor=south]{$s_{AB}$} (process15);
\draw [arrow1] (process5) -- node[anchor=east] {$r_{10}$} (process1);
\draw [arrow1] (process6) -- node[anchor=east] {$r_{11}$} (process2);
\draw [arrow1] (process7) -- node[anchor=east] {$r_{1\cdot}$} (process3);
\draw [arrow1] (process8) -- node[anchor=east] {$r_{1B}$} (process4);
\draw [arrow1] (process9) -- node[anchor=east] {$r_{20}$} (process5);
\draw [arrow1] (process10) -- node[anchor=east] {$r_{21}$} (process6);
\draw [arrow1] (process11) -- node[anchor=east] {$r_{2\cdot}$} (process7);
\draw [arrow1] (process12) -- node[anchor=east] {$r_{2B}$} (process8);
\draw [arrow1] (process13) -- node[anchor=east] {$r_{A1}$} (process9);
\draw [arrow1] (process14) -- node[anchor=east] {$r_{A2}$} (process10);
\draw [arrow1] (process15) -- node[anchor=east] {$r_{A\cdot}$} (process11);
\draw [arrow1] (process16) -- node[anchor=east] {$r_{AB}$} (process12);
\end{tikzpicture}
\end{CJK*}
\eg

Once we desire to know $F(A,B)=(?)F(0,0)+\mbox {etc.}$, this graph will be helpful.

The number next to an arrow is called a ``weight'', meaning the factor provided by the previous term to the next term. A propagation occurs only along the arrow. \par

We define the ``weight'' from $F(A,B)$ to $F(0,0)$  in the following manner:
\begin{itemize}

\item[(i)] When $F(A,B)=F(0,0)$, we define the weight to be $1$.

\item[(ii)] When there is no viable path from $F(A,B)$ to $F(0,0)$  and $F(A,B)\neq F(0,0)$, we define the weight to be $0$.

  \item[(iii)] When there are not the same and there exists at least one viable path from $F(A,B)$ to $F(0,0)$, we first find all such paths, then for every such path, multiply all weights along it to get a weight (which is called the path weight for this path), and finally, add up path weights for all paths mentioned above to get a new weight, which is the weight from
  $F(A,B)$ to $F(0,0)$.

\end{itemize}

Based on Graph \ref{graph1}, in this paper we will introduce some back-propagation graphs (BPG(s) for short), in the sense of integration by parts, to give the proof of Carleman estimate. Through a BPG, one can easily analyse the weight relationship, allowing one obtains the results of integration by parts without any practical calculation.

\section{BPM for the estimation of $I_1I_2$}\label{sec3.1}

We arbitrarily fix an integer $n\ge2$. For any integer $m$ with $0\le m\le n-1$, we hope to find the coefficient of the final term $\ell_x^{2n-2m-2}\ell_{xx}|\pa_x^mw|^2$ in $I_1I_2$. Set
$$p\=2n-2m-1.$$

We introduce the following back-propagation graph:
\bg\label{graph5}
\begin{CJK*}{GBK}{kai}
\ \\
\tikzstyle{process} = [rectangle,rounded corners, minimum width=1.1cm,minimum height=0.8cm,text centered,text width =2.9cm, draw=black,fill=white]
\tikzstyle{arrow1} = [thick,->,>=stealth]
{$\\ \mbox{}$}
\begin{tikzpicture}
\node (process1) [process] {$\ell_x^{p-1}\ell_{xx}\pa_x^mw\pa_x^mw$};
\node (process2) [process,right of=process1,xshift=3.4cm,text width =3.6cm] {$\ell_x^{p-1}\ell_{xx}\pa_x^{m+1}w\pa_x^{m-1}w$};
\node (process3) [process,right of=process2,xshift=2.3cm,text width =0.8cm] {...};
\node (process4) [process,right of=process3,xshift=2.0cm] {$\ell_x^{p-1}\ell_{xx}\pa_x^{2m}w\pa_x^0w$};
\node (process5) [process,below of=process1,yshift=-0.3cm]{$\ell_x^{p}\pa_x^{m+1}w\pa_x^mw$};
\node (process6) [process,right of=process5,xshift=3.4cm,text width =3.6cm] {$\ell_x^{p}\pa_x^{m+2}w\pa_x^{m-1}w$};
\node (process7) [process,right of=process6,xshift=2.3cm,text width =0.8cm] {...};
\node (process8) [process,right of=process7,xshift=2.0cm] {$\ell_x^{p}\pa_x^{2m+1}w\pa_x^0w$};
\draw [arrow1] (process2) -- node[anchor=south]{$-1$} (process1);
\draw [arrow1] (process3) -- node[anchor=south]{$-1$} (process2);
\draw [arrow1] (process4) -- node[anchor=south]{$-1$} (process3);
\draw [arrow1] (process6) -- node[anchor=south]{$-1$} (process5);
\draw [arrow1] (process7) -- node[anchor=south]{$-1$} (process6);
\draw [arrow1] (process8) -- node[anchor=south]{$-1$} (process7);
\draw [arrow1] (process5) -- node[anchor=east] {$-p/2$} (process1);
\draw [arrow1] (process6) -- node[anchor=east] {$-p$} (process2);
\draw [arrow1] (process7) -- node[anchor=east] { } (process3);
\draw [arrow1] (process8) -- node[anchor=east] {$-p$} (process4);
\end{tikzpicture}
\end{CJK*}
\eg

The reason that the above graph is called the ``back-propagation graph''  is that we start drawing this graph from the goal term which is asked for, and the rest part of the graph is drawn against the direction of the propagation.

Let $A$ and $B$ be two terms in two adjacent positions. An arrow from $A$ to $B$ with the number $k$ means that  $A=kB+etc.$, through one step of integration by parts. For example, ``$\ell_x^{p-1}\ell_{xx}\pa_x^{m+1}w\pa_x^{m-1}w\stackrel{-1}\longrightarrow\ell_x^{p-1}\ell_{xx}\pa_x^mw\pa_x^mw$'' means
\be\ba{ll}
\ds\ell_x^{p-1}\ell_{xx}\pa_x^{m+1}w\pa_x^{m-1}w \\
\ns\ds=(\ell_x^{p-1}\ell_{xx}\pa_x^mw\pa_x^mw)_x-\ell_x^{p-1}\ell_{xx}\pa_x^mw\pa_x^mw-(p-1)\ell_x^{p-2}\ell_{xx}^2\pa_x^mw\pa_x^{m-1}w\\
\ns\ds=(-1)\cdot\ell_x^{p-1}\ell_{xx}\pa_x^mw\pa_x^mw+etc..
\ea\eee

Based on Graph \ref{graph5}, we have the following result.

 \bp\label{22051602}
Under the assumption of Proposition \ref{lmSW02}, we have
\be
\ds I_1I_2=(\cdot)_x+\sum_{m=0}^{n-1}\[\frac{n^2}{2}C_{n-1}^{m}\ell_x^{2n-2m-2}\ell_{xx}+O(\l^{2n-2m-3})\]|\pa_x^{m}w|^2.
\eee
\ep

\br
Based on the selection of $I_1,I_2$, one knows that $I_1I_2$ does not contain such terms $\ell_x^{2n-2m-1-s}\ell_{xx}^s|\pa_x^{m}w|^2$ with $s$ even, after integration by parts.
\er

Before giving the proof of Proposition \ref{22051602},  we first give an example to show how the BPG works for calculating the coefficient of some term in $I_1I_2$ through integration by parts.

\begin{example}\label{exam1} For $n=4$, find the coefficient of $\ell_x^2\ell_{xx}\pa_x^2w\pa_x^2w$ in $I_1I_2$ after integration by parts. By choosing $p=3, m=2$ in Graph \ref{graph5}, we have the following BPG:
\bg\label{graph3}
\begin{CJK*}{GBK}{kai}\fontsize{10.5}{6}\selectfont
\ \\
\tikzstyle{process} = [rectangle,rounded corners, minimum width=2cm,minimum height=0.7cm,text centered,text width =3.7cm, draw=black,fill=white]
\tikzstyle{arrow1} = [thick,->,>=stealth]
{$\\ \mbox{}$}
\begin{tikzpicture}
\node (process1) [process] {$\ell_x^2\ell_{xx}\pa_x^2w\pa_x^2w$};
\node (process2) [process,right of=process1,xshift=3.7cm] {$\ell_x^2\ell_{xx}\pa_x^3w\pa_x^1w$};
\node (process3) [process,right of=process2,xshift=3.7cm] {$\ell_x^2\ell_{xx}\pa_x^4w\pa_x^0w$};
\node (process5) [process,below of=process1,yshift=-.3cm]{$\ell_x^3\pa_x^3w\pa_x^2w$};
\node (process6) [process,right of=process5,xshift=3.7cm] {$\ell_x^3\pa_x^4w\pa_x^1w$};
\node (process7) [process,right of=process6,xshift=3.7cm] {$\ell_x^3\pa_x^5w\pa_x^0w$};
\draw [arrow1] (process2) -- node[anchor=south]{$-1$} (process1);
\draw [arrow1] (process3) -- node[anchor=south]{$-1$} (process2);
\draw [arrow1] (process6) -- node[anchor=south]{$-1$} (process5);
\draw [arrow1] (process7) -- node[anchor=south]{$-1$} (process6);
\draw [arrow1] (process5) -- node[anchor=east] {$-3/2$} (process1);
\draw [arrow1] (process6) -- node[anchor=east] {$-3$} (process2);
\draw [arrow1] (process7) -- node[anchor=east] {$-3$} (process3);
\end{tikzpicture}
\end{CJK*}
\eg
\ \\
Then fill the coefficients in $I_1I_2$ from Example \ref{I1I2n4}, we have
\bg\label{graph4}
\begin{CJK*}{GBK}{kai}\fontsize{10.5}{6}\selectfont
\ \\
\tikzstyle{process} = [rectangle,rounded corners, minimum width=2cm,minimum height=0.7cm,text centered,text width =3.7cm, draw=black,fill=white]
\tikzstyle{arrow1} = [thick,->,>=stealth]
{$\\ \mbox{}$}
\begin{tikzpicture}
\node (process1) [process] {$(-6\times 6)\ell_x^2\ell_{xx}\pa_x^2w\pa_x^2w$};
\node (process2) [process,right of=process1,xshift=3.7cm] {$(-4\times 12)\ell_x^2\ell_{xx}\pa_x^3w\pa_x^1w$};
\node (process3) [process,right of=process2,xshift=3.7cm] {$(-6\times 1)\ell_x^2\ell_{xx}\pa_x^4w\pa_x^0w$};
\node (process5) [process,below of=process1,yshift=-.3cm]{$(-4\times 6)\ell_x^3\pa_x^3w\pa_x^2w$};
\node (process6) [process,right of=process5,xshift=3.7cm] {$(-4\times 1)\ell_x^3\pa_x^4w\pa_x^1w$};
\node (process7) [process,right of=process6,xshift=3.7cm] {$0\cdot \ell_x^3\pa_x^5w\pa_x^0w$};
\draw [arrow1] (process2) -- node[anchor=south]{$-1$} (process1);
\draw [arrow1] (process3) -- node[anchor=south]{$-1$} (process2);
\draw [arrow1] (process6) -- node[anchor=south]{$-1$} (process5);
\draw [arrow1] (process7) -- node[anchor=south]{$-1$} (process6);
\draw [arrow1] (process5) -- node[anchor=east] {$-3/2$} (process1);
\draw [arrow1] (process6) -- node[anchor=east] {$-3$} (process2);
\draw [arrow1] (process7) -- node[anchor=east] {$-3$} (process3);
\end{tikzpicture}
\end{CJK*}
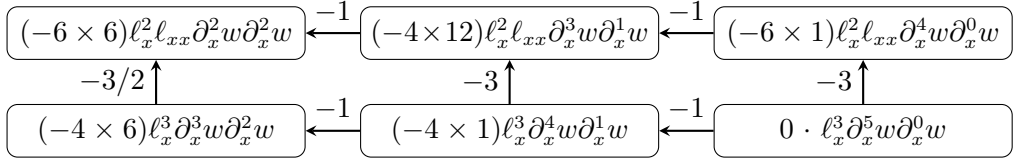
\eg
\ \\

Though there are many terms in $I_1I_2$, only at most $6$ terms contributing to $\ell_x^2\ell_{xx}\pa_x^2w\pa_x^2w$. Thus we only need to see how these $6$ terms contribute coefficients to
$\ell_x^2\ell_{xx}\pa_x^2w\pa_x^2w$.

We write $\ds d(r, s, a, b)=\sum_{D} h(r,s,a,b)g(r,s,a,b)$ by the coefficient of $\ell_x^r\ell_{xx}^s\pa_x^aw\pa_x^bw$, where $D=\{(r,s,a,b)\in\mathbb{Z}_{\ge0}^4|r+2s+a+b=8,\ s=0,1,\ a\ge b\}$, $h(r,s,a,b)$ is the coefficient of each term, and $g(r,s,a,b)$ is the weight from each term to the goal term. One can know $g(r,s,a,b)$ from the BPG and $h(r,s,a,b)$ from $I_1,I_2$ in (\ref{22053001}).

Here, we try to compute $d(2,1,2,2)$, i.e., the final coefficient of $\ell_x^2\ell_{xx}\pa_x^2w\pa_x^2w$ in $I_1I_2$. By the BPG, we know that
 $$\ba{ll}\ds D=\{(r,s,a,b)\}=\big\{(2,1,2,2),(2,1,3,1),(2,1,4,0),\\
\ns\ds\qq\qq\q\qq\qq(3,0,3,2),(3,0,4,1),(3,0,5,0)\big\},\ea$$
\be\left\{\begin{matrix}
\ns &h(2,1,2,2)=&-6\times6,  &h(2,1,3,1)=&-4\times12, &h(2,1,4,0)=&-6\times1,\\
\ns &h(3,0,3,2)=&-4\times6, &h(3,0,4,1)=&-4\times1, &h(3,0,5,0)=& 0,
\end{matrix}\right.\eee
and
\be\left\{\begin{matrix}
\ns &g(2,1,2,2)=&1,\\
\ns &g(2,1,3,1)=&(-1)^{1},\\
\ns &g(2,1,4,0)=&(-1)^{2},\\
\ns &g(3,0,3,2)=&\ds(-1)^{1}3/2,\\
\ns &g(3,0,4,1)=&\ds(-1)^{2}(3/2+3),\\
\ns &g(3,0,5,0)=&\ds(-1)^{3}(3/2+3+3).
\end{matrix}\right.\eee
Then:
\be\ba{ll}
d(2,1,2,2)&\ds=(-36)\cdot 1+(-48)(-1)^1+(-6)(-1)^2 \\
&\q+(-24)(3/2)(-1)^1+(-4)(3/2+3)(-1)^2+0\cdot(3/2+3+3)(-1)^3\\
&=24,
\ea\eee
which means
\be\begin{matrix}
-6\times 6&\ell_x^2\ell_{xx}\pa_x^2w\pa_x^2w&= &-36&\ell_x^2\ell_{xx}|\pa_x^2w|^2+etc., \\
-4\times 12&\ell_x^2\ell_{xx}\pa_x^3w\pa_x^1w&= &48&\ell_x^2\ell_{xx}|\pa_x^2w|^2+etc., \\
-6\times 1&\ell_x^2\ell_{xx}\pa_x^4w\pa_x^0w&= &-6&\ell_x^2\ell_{xx}|\pa_x^2w|^2+etc., \\
-4\times 6&\ell_x^3\pa_x^3w\pa_x^2w&= &36&\ell_x^2\ell_{xx}|\pa_x^2w|^2+etc., \\
-4\times 1&\ell_x^3\pa_x^4w\pa_x^1w&= &-18&\ell_x^2\ell_{xx}|\pa_x^2w|^2+etc., \\
0&\ell_x^3\pa_x^5w\pa_x^0w&= &0&\ell_x^2\ell_{xx}|\pa_x^2w|^2+etc., \\
\Rightarrow &\ Total&= &24&\ell_x^2\ell_{xx}|\pa_x^2w|^2+etc.. \\
\end{matrix}\eee
\end{example}

\ \\
{\bf Proof of Proposition \ref{22051602}. } Example \ref{exam1} (when $n=4$) has explained our main idea. Here we give the detailed proof for an arbitrarily fixed  $n$. We write $d_m$ (with $0\le m\le n-1$) for the coefficient of the term $\ell_x^{p-1}\ell_{xx}|\pa_x^m w|^2$ for $0\le m\le n-1$. Then we have
\bel{22051609}
\ds d_m=\sum_D h(r,s,a,b)g(r,s,a,b),
\ee
where $D=\{(r,s,a,b)\in\mathbb{Z}_{\ge0}^4|r+2s+a+b=2n,\ s=0,1,\ a\ge b\}$, $h(r,s,a,b)$ is the coefficient of $\ell_x^{r}\ell_{xx}^s\pa_x^aw\pa_x^bw$ in $I_1I_2$ before integration by parts, and $g(r,s,a,b)$ is the weight from $\ell_x^{r}\ell_{xx}^s\pa_x^aw\pa_x^bw$ to $\ell_x^{p-1}\ell_{xx}|\pa_x^mw|^2$ in Graph \ref{graph5}. \par
On the one hand, by (\ref{1128-b0}), comparing the terms in $I_1I_2$ to find their coefficients, one can get
\bel{22051610}
h(r,s,a,b)=\left\{\ba{ll}
\ns\ds (-1)^{a+b}(-C_{n}^{2}C_{n}^{m}C_{n-2}^{m}),\ s=1,\ a=b=m, \\
\ns\ds (-1)^{a+b}[-C_{n}^{2}(C_{n}^{a}C_{n-2}^{b}+C_{n}^{b}C_{n-2}^{a})],\\
\ns\ds\qq\ s=1,\ m+1\le a\le 2m,\ b=2m-a; \\
\ns\ds (-1)^{a+b}(C_{n}^{a}C_{n}^{b}),\\
\ns\ds\qq\ s=0,\ m+1\le a\le 2m+1,\ b=2m+1-a.
\ea\right.\ee \par
On the other hand, from the BPG one can obtain the weight:
\bel{22051611}
g(r,s,a,b)=\left\{\ba{ll}
\ns\ds (-1)^{a-m},\ s=1,\ m\le a\le 2m,\\
\ns\ds (-1)^{a-m}\frac{p}{2}[2(a-m)-1],\\
\ns\ds\qq\qq\q s=0,\ m+1\le a\le 2m+1.
\ea\right.\ee \par
Combining (\ref{22051609})-(\ref{22051611}), one has
\bel{22051612}\ba{ll}
d_m&\ds=(-1)^{m-m}(-1)^{m+m}(-C_{n}^{2}C_{n}^{m}C_{n-2}^{m}) \\
&\ds\q +\sum_{a=m+1}^{2m}(-1)^{a-m}(-1)^{a+2m-a}[-C_{n}^{2}(C_{n}^{a}C_{n-2}^{2m-a}+C_{n}^{2m-a}C_{n-2}^{a})] \\
&\ds\q +\sum_{a=m+1}^{2m+1}(-1)^{a-m}\frac{p}{2}[2(a-m)-1](-1)^{a+2m+1-a}(C_{n}^{a}C_{n}^{2m+1-a}) \\
&\ds=-C_{n}^{2}C_{n}^{m}C_{n-2}^{m}-C_{n}^{2}\sum_{k=1}^{m}(-1)^{k}(C_{n}^{m+k}C_{n-2}^{m-k}+C_{n}^{m-k}C_{n-2}^{m+k}) \\
&\ds\q +\frac{p}{2}\sum_{k=0}^{m}(-1)^{k}(1+2k)C_{n}^{m+1+k}C_{n}^{m-k}.
\ea\ee
Then, by Proposition \ref{22051601}, we are led to that $\ds d_m=\cK_{n,m}=\frac{n^2}{2}C_{n-1}^m$.
\endpf

\section{BPM for the estimation of $w_tI_2$ }

First, by (\ref{1128-b0}) and through integration by parts, we have
\bel{1207-a0}\ba{ll}
\ns w_tI_2=&\ds(\cdot)_t+(\cdot)_x+\sum_{r+2s+2m=n}(\cdot)(\ell_x^r\ell_{xx}^s)_t|\pa_x^mw|^2 \\
\ns&\ds +\sum_{D}d_C(r,s,m)\ell_x^{r}\ell_{xx}^{s}\pa_x^{m}\pa_t w\pa_x^{m+1}w,
\ea\ee
where\bel{24012101}
\ds D=\{(r,s,m)\in\mathbb{Z}_{\ge0}^3|r+2s+2m+1=n\}.
\ee
Since $\ell_{tx}=0$, we have $(\ell_x^r\ell_{xx}^s)_t=0$. Therefore, we only need to calculate the coefficient $d_C$.

To calculate $d_C$, we introduce the following BPG:
\bg\label{graph6}
\begin{CJK*}{GBK}{kai}\fontsize{9}{6}\selectfont
\ \\
\tikzstyle{process} = [rectangle,rounded corners, minimum height=0.8cm,text centered, text width=3.3cm, draw=black,fill=white]
\tikzstyle{arrow1} = [thick,->,>=stealth]
{$\\ \mbox{}$}
\begin{tikzpicture}
\node (process1) [process]{$\ell_x^r\ell_{xx}^s\pa_x^{m}\pa_t w\pa_x^{m+1}w$};
\node (process2) [process,right of=process1,xshift=3.6cm,text width=4.0cm] {$\ell_x^r\ell_{xx}^s\pa_x^{m-1}\pa_t w\pa_x^{m+2}w$};
\node (process3) [process,right of=process2,xshift=2.2cm,text width=0.6cm] {...};
\node (process4) [process,right of=process3,xshift=1.9cm] {$\ell_x^r\ell_{xx}^s\pa_t w\pa_x^{2m+1}w$};
\node (process5) [process,below of=process1,yshift=-.4cm]{$\ell_x^{r+1}\ell_{xx}^{s-1}\pa_x^{m}\pa_t w\pa_x^{m+2}w$};
\node (process6) [process,right of=process5,xshift=3.6cm,text width=4.0cm] {$\ell_x^{r+1}\ell_{xx}^{s-1}\pa_x^{m-1}\pa_t w\pa_x^{m+3}w$};
\node (process7) [process,right of=process6,xshift=2.2cm,text width=0.6cm] {...};
\node (process8) [process,right of=process7,xshift=1.9cm] {$\ell_x^{r+1}\ell_{xx}^{s-1}\pa_t w\pa_x^{2m+2}w$};
\node (process9) [process,below of=process5,yshift=-.4cm]{...};
\node (process10) [process,right of=process9,xshift=3.6cm,text width=4.0cm] {...};
\node (process11) [process,right of=process10,xshift=2.2cm,text width=0.6cm] {...};
\node (process12) [process,right of=process11,xshift=1.9cm] {...};
\node (process13) [process,below of=process9,yshift=-.4cm]{$\ell_x^{r+s-1}\ell_{xx}^{1}\pa_x^{m}\pa_t w\pa_x^{m+s}w$};
\node (process14) [process,right of=process13,xshift=3.6cm,text width=4.0cm] {$\ell_x^{r+s-1}\ell_{xx}^{1}\pa_x^{m-1}\pa_t w\pa_x^{m+s+1}w$};
\node (process15) [process,right of=process14,xshift=2.2cm,text width=0.6cm] {...};
\node (process16) [process,right of=process15,xshift=1.9cm] {$\ell_x^{r+s-1}\ell_{xx}^{1}\pa_t w\pa_x^{2m+s}w$};
\node (process17) [process,below of=process13,yshift=-.4cm]{$\ell_x^{r+s}\pa_x^{m}\pa_t w\pa_x^{m+s+1}w$};
\node (process18) [process,right of=process17,xshift=3.6cm,text width=4.0cm] {$\ell_x^{r+s}\pa_x^{m-1}\pa_t w\pa_x^{m+s+2}w$};
\node (process19) [process,right of=process18,xshift=2.2cm,text width=0.6cm] {...};
\node (process20) [process,right of=process19,xshift=1.9cm] {$\ell_x^{r+s}\pa_t w\pa_x^{2m+s+1}w$};
\draw [arrow1] (process2) -- node[anchor=south]{$-1$} (process1);
\draw [arrow1] (process3) -- node[anchor=south]{$-1$} (process2);
\draw [arrow1] (process4) -- node[anchor=south]{$-1$} (process3);
\draw [arrow1] (process6) -- node[anchor=south]{$-1$} (process5);
\draw [arrow1] (process7) -- node[anchor=south]{$-1$} (process6);
\draw [arrow1] (process8) -- node[anchor=south]{$-1$} (process7);
\draw [arrow1] (process10) -- node[anchor=south]{$-1$} (process9);
\draw [arrow1] (process11) -- node[anchor=south]{$-1$} (process10);
\draw [arrow1] (process12) -- node[anchor=south]{$-1$} (process11);
\draw [arrow1] (process14) -- node[anchor=south]{$-1$} (process13);
\draw [arrow1] (process15) -- node[anchor=south]{$-1$} (process14);
\draw [arrow1] (process16) -- node[anchor=south]{$-1$} (process15);
\draw [arrow1] (process18) -- node[anchor=south]{$-1$} (process17);
\draw [arrow1] (process19) -- node[anchor=south]{$-1$} (process18);
\draw [arrow1] (process20) -- node[anchor=south]{$-1$} (process19);
\draw [arrow1] (process5) -- node[anchor=east] {$-(r+1)$} (process1);
\draw [arrow1] (process6) -- node[anchor=east] {$-(r+1)$} (process2);
\draw [arrow1] (process7) -- node[anchor=east] {$-(r+1)$} (process3);
\draw [arrow1] (process8) -- node[anchor=east] {$-(r+1)$} (process4);
\draw [arrow1] (process9) -- node[anchor=east] {$-(r+2)$} (process5);
\draw [arrow1] (process10) -- node[anchor=east] {$-(r+2)$} (process6);
\draw [arrow1] (process11) -- node[anchor=east] {$-(r+2)$} (process7);
\draw [arrow1] (process12) -- node[anchor=east] {$-(r+2)$} (process8);
\draw [arrow1] (process13) -- node[anchor=east] {$-(r+s-1)$} (process9);
\draw [arrow1] (process14) -- node[anchor=east] {$-(r+s-1)$} (process10);
\draw [arrow1] (process15) -- node[anchor=east] {$-(r+s-1)$} (process11);
\draw [arrow1] (process16) -- node[anchor=east] {$-(r+s-1)$} (process12);
\draw [arrow1] (process17) -- node[anchor=east] {$-(r+s)$} (process13);
\draw [arrow1] (process18) -- node[anchor=east] {$-(r+s)$} (process14);
\draw [arrow1] (process19) -- node[anchor=east] {$-(r+s)$} (process15);
\draw [arrow1] (process20) -- node[anchor=east] {$-(r+s)$} (process16);
\end{tikzpicture}
\end{CJK*}
\eg

Based on Graph \ref{graph6}, we have the following result.

\bp\label{23052501}
Under the assumption of Proposition \ref{lmSW02},
\bel{0123-a}
w_tI_2(w)=(\cdot)_t+(\cdot)_x+\sum_{D_C}d_C(r,s,m)\ell_x^{r}\ell_{xx}^{s}\pa_x^{m}\pa_t w\pa_x^{m+1}w,
\ee
where $\ds D_C\=\{(r,s,m)\in\mathbb{Z}_{\ge0}^3|r+2s+2m+1=n,\ s\ge3, (-1)^s=-1\}$ and
$$\ds d_C\=(-1)^{n+m}\frac{(s-1)(2m+s)(n-2m-s-1)}{2(m+s)}C_n^{2m+s+1}C_{m+s}^{s}\prod_{j=1}^{s-1}(r+j).$$
\ep
\ \\
\br In the case that $n\le 6$, we know that $D_C=\varnothing$. In the case that $n=7$, we know that $D_C=\{(0,3,0)\}$ and  $w_tI_2(w)=(\cdot)_t+(\cdot)_x-210\ell_{xx}^3\pa_t w\pa_x w$.
\er
\ \\
{\bf Proof of Proposition \ref{23052501}.}  Similar to what we did in the proof of Proposition \ref{22051602}, we write
$$\ds d_C=\sum_{D'}h(r,s,a,b)g(r,s,a,b),$$
where the domain $D'$ is the set of all parameters sets $(r,s,a,b)$ by terms $\ell_x^r\ell_{xx}^s\pa_x^{a}\pa_t w\pa_x^{b}w$ in Graph \ref{graph6}, and $h(r,s,a,b)$ and $g(r,s,a,b)$ denote the coefficients and weights of term $\ell_x^r\ell_{xx}^s\pa_x^a\pa_tw\pa_x^bw$, respectively. Recall that $I_2(w)$ has the following form:
$$\ba{ll}\ds
I_2(w)=\sum_{even\ k\in[0,n]}C_{n}^{k}(-1)^{n-k}\ell_x^{n-k}\pa_x^{k}w \\
\ns\ds\q\q\qq -C_{n}^2\ell_{xx}\sum_{odd\ j\in[0,n-2]}C_{n-2}^{j}(-1)^{n-2-j}\ell_x^{n-2-j}\pa_x^{j}w.
\ea$$
Hence, we can find at most two terms $\ell_x^{r+s-1}\ell_{xx}^{1}\pa_t w\pa_x^{2m+s}w$ and $\ell_x^{r+s}\pa_t w\pa_x^{2m+s+1}w$, with possible non-zero coefficients, in the bottom right corner of Graph \ref{graph6}. Then we see that
$$\ba{ll}\ds d_C(r,s,m)=&\ds h(r+s-1,1,0,2m+s)g(r+s-1,1,0,2m+s)\\
\ns&\ds+h(r+s,0,0,2m+s+1)g(r+s,0,0,2m+s+1).\ea$$
If $s$ is even, from $I_2$ and the BPG we know
\bel{1207-0}
h(r+s-1,1,0,2m+s)=h(r+s,0,0,2m+s+1)=0.
\ee
If $s$ is odd, we can find
\bel{23092108}\left\{\ba{ll}
\ns\ds h(r+s-1,1,0,2m+s)=(-1)^n C_n^2 C_{n-2}^{2m+s}, \\
\ns\ds g(r+s-1,1,0,2m+s)=(-1)^{m+(s-1)} C_{m+(s-1)}^m \prod_{j=1}^{s-1}(r+j), \\
\ns\ds h(r+s,0,0,2m+s+1)=(-1)^n C_n^{2m+s+1}, \\
\ns\ds g(r+s,0,0,2m+s+1)=(-1)^{m+s} C_{m+s}^m \prod_{j=1}^{s}(r+j),
\ea\right.\ee
therefore
\bel{0123-b}\ba{ll}
\ns &\ds d_C(r,s,m) \\
\ns &\ds =(-1)^{n+m}[C_{s+m-1}^{s-1}C_{n}^{2}C_{n-2}^{2m+s}-C_{s+m}^{s}(r+s)C_{n}^{2m+s+1}]\prod_{j=1}^{s-1}(r+j) \\
\ns &\ds =(-1)^{n+m}\frac{(s-1)(2m+s)(n-2m-s-1)}{2(m+s)}C_n^{2m+s+1}C_{m+s}^{s}\prod_{j=1}^{s-1}(r+j).
\ea\ee
In the last equation, we used the fact that $r=n-2s-2m-1$, which follows from $(r,s,m)\in D$ defined in (\ref{24012101}). Set $D_C\=D\cap\{s\ge3,\ (-1)^s=-1\}$. By (\ref{1207-0}) and (\ref{0123-b}), we know that $$
\ds \sum_{D}d_C(r,s,m)\ell_x^{r}\ell_{xx}^{s}\pa_x^{m}\pa_t w\pa_x^{m+1}w=\sum_{D_C}d_C(r,s,m)\ell_x^{r}\ell_{xx}^{s}\pa_x^{m}\pa_t w\pa_x^{m+1}w.$$
This, along with (\ref{1207-a0}), implies (\ref{0123-a}). \endpf

\section{ Proof of Theorem \ref{thmSW00}}

This section proves Theorem \ref{thmSW00}.

\ \\
{\bf Proof. } By the definitions  of $I_1$ and $I_2$ (see (\ref{1128-b0})), using Propositions \ref{22051602}, \ref{23052501}, integrating (\ref{23061404}) over $(0,T)\t(0,L)$, we obtain that when $v\in C_0^\i([0,T]\t[0,L])$ and $w=\th v$,
\bel{1208}\ba{ll}\ds
\int_0^T\int_0^L \th^2|\cP v|^2dxdt\\
\ns\ds\ge \int_0^T\int_0^L\[\sum_{m=0}^{n-1}n^2C_{n-1}^m\ell_x^{2n-2m-2}\ell_{xx}+O(\l^{2n-2m-3})\]|\pa_x^mw|^2dxdt\\
\ns\ds\q+2\int_0^T\int_0^L \a \sum_{D_C}d_C(r,s,m)\ell_x^{r}\ell_{xx}^{s}\pa_x^{m}\pa_t w\pa_x^{m+1}w dxdt,
\ea\ee
where the last term $\ds\sum_{D_C}d_C(r,s,m)\ell_x^{r}\ell_{xx}^{s}\pa_x^{m}\pa_t w\pa_x^{m+1}w$ is defined in Proposition \ref{23052501}. Noting that
\bel{23120201}\ba{ll}
\ns\ds\ell_x^r\ell_{xx}^s\pa_x^a\pa_tw\pa_x^{m+1}w=&\ds(\ell_x^r\ell_{xx}^s\pa_x^{a-1}\pa_tw\pa_x^{m+1}w)_x-(\ell_x^r\ell_{xx}^s)_x\pa_x^{a-1}\pa_tw\pa_x^{m+1}w\\
\ns&\ds-\ell_x^r\ell_{xx}^s\pa_x^{a-1}\pa_tw\pa_x^{m+2}w,
\ea\ee
with $(r,s,a,m)\in D_1\=\{(r,s,a,m)\in\mathbb{Z}_{\ge0}^4|1\le a\le m,\ r+2s+a+m+1=n,\ s\ge3\}$, where $a$ decreases to $1$ during the iteration, one knows
\be
\ds\ell_x^r\ell_{xx}^s\pa_x^a\pa_tw\pa_x^{m+1}w=(\cdot)_x+\pa_tw\sum_{D_2}(\cdot)\ell_x^r\ell_{xx}^s\pa_x^mw,
\eee
where $D_2\=\{(r,s,m)\in\mathbb{Z}_{\ge0}^3|r+2s+m=n,\ s\ge3\}$.
Therefore,
\bel{23120202}\ba{ll}
\ns\ds \sum_{D_C}(\cdot)\ell_x^{r}\ell_{xx}^{s}\pa_x^{m}\pa_tw\pa_x^{m+1}w &\ds=(\cdot)_x+\sum_{D_1}(\cdot)\ell_x^r\ell_{xx}^s\pa_x^a\pa_tw\pa_x^{m+1}w \\
\ns&\ds=(\cdot)_x+\pa_tw\sum_{D_2}(\cdot)\ell_x^r\ell_{xx}^s\pa_x^mw.
\ea\ee

For convenience, we set $\ds\Phi\=\sum_{D_2}(\cdot)\ell_x^r\ell_{xx}^s\pa_x^mw$ in (\ref{23120202}). Then by (\ref{as00}), we find
\bel{1208-a2}\ba{ll}\ds
2\a w_t\Phi&\ds=2\Phi\[\th\cP v+\a\ell_t w-I_1-I_2-I_3\]\\
\ns&\ds\ge-\th^2|\cP v|^2-|\Phi|^2+2\Phi\[\a\ell_t w-I_1-I_2-I_3\].
\ea\ee

In the case that $s\ge 3$, $\Phi$ is a lower order term compared with $I_1$ or $I_2$. Comparing the order of $\ell_{xx}$, we can obtain
\bel{1208-c0}
\Phi (\a\ell_t w+I_1+I_2+I_3)=O(\l^{2n-2m-3})|\pa_x^mw|^2.
\ee
Now, it follows from (\ref{1208}) and (\ref{1208-c0}) that
$$\ba{ll}\ds
\int_0^T\int_0^L \th^2|\cP v|^2dxdt\\
\ns\ds\ge \int_0^T\int_0^L\[\sum_{m=0}^{n-1}n^2C_{n-1}^m\ell_x^{2n-2m-2}\ell_{xx}+O(\l^{2n-2m-3})\]|\pa_x^mw|^2dxdt.
\ea$$
This completes the proof of Theorem \ref{thmSW00}.
\endpf

\section{Application}\label{section4}

In this section, we will prove the  conditional stability  in a Cauchy problem  for a
time fractional diffusion equation with $\frac{1}{3}$-order with the help of Theorem \ref{thmSW00}.
First, we transfer the original equation to an usual p.d.e. with its principal part $\pa_t-\pa_x^6$, with the aid of the Caputo derivative. Then we apply our Theorem \ref{thmSW00} to obtain a Carleman-type inequality for the new equation. Finally, we use the similar way to that used in \cite{XCY} to obtain the desired conditional stability.

Given $T>0$, $L>0$, we consider the following fractional diffusion equation:
\bel{23092101}\left\{\ba{ll}
\ns\ds \pa_t^{1/3}u(t,x)-\pa_x^2u(t,x)=f(t,x),&\ (t,x)\in (0,T)\times(0,L),\\
\ns\ds u(0,x)=0,&\ x\in(0,L), \\
\ns\ds \pa_x^ju(t,x)|_{x=0}=h_j(t) ( j=0, 1), &\ t\in (0,T),
\ea\right.\ee
where $f(x)$ and $h_j(t)$ ($j=0, 1)$ are given, $\pa_t^{1/3}$ is the fractional order derivative in the Caputo sense which is defined by
$$\ds \pa_t^{\g}y(t)\=\frac{1}{\G(1-\g)}\int_0^t\frac{\pa_sy(s)}{(t-s)^\g}ds,\ \g\in(0,1),$$
and $\ds\Gamma(x)=\int_0^{+\i}s^{x-1}e^{-s}ds$ is the Gamma Function.

Before the main result, we introduce some notation. Given  $0<\delta_0<T/2$, we let
$$Q\=\{(t,x)|\delta_0<t<T,\ 0<x<L\},\ \O\=(0,L).$$
Given $0<\e<\min\{L^2,1\}$ and  $L<x_0<L+\sqrt{\e}/2$, we define
$$Q_{\e}\=\{(t,x)\in Q|\psi(t,x)>\e\},$$
where
$$\psi\=(x-x_0)^2-\beta(t-T/2)^2$$
with $\beta>x_0^2/(\delta_0-T/2)^2$. One can easily see that  $Q_{\e}$ is a nonempty subset of $Q$ when $\e$ is small.  By (\ref{23092101}), one can easily check that if  $u\in C^{1,6}(\cl{Q})$, then
\bel{23101110}\left\{\ba{ll}
 \pa_x^2u(t,0)=\pa_t^{1/3}h_0(t)-f(t,0)\=h_2(t),\\
\ns\ds \pa_x^3u(t,0)=\pa_t^{1/3}h_1(t)-\pa_xf(t,0)\=h_3(t),\\
\ns\ds \pa_x^4u(t,0)=\pa_t^{1/3}\pa_t^{1/3}h_0(t)-\pa_t^{1/3}f(t,0)-\pa_x^2f(t,0)\=h_4(t),\\
\ns\ds \pa_x^5u(t,0)=\pa_t^{1/3}\pa_t^{1/3}h_1(t)-\pa_t^{1/3}\pa_xf(t,0)-\pa_x^3f(t,0)\=h_5(t).
\ea\right.\ee

\bt\label{thmapp}
Assume that $u\in C^{1,6}(\cl{Q})$ satisfies (\ref{23092101}), and $f\in C^1([0,T];L^2(\O))\cap L^2(0,T;H^{4}(\O))$. Let $\psi(t,x)$ be given chosen such that $Q_{\e/4}\subset Q$ and $h_j (j=2, 3, 4, 5)$ be given by (\ref{23101110}). Then there exist $C>0$ and $\tau\in(0,1)$ such that
\bel{23101001}
\sum_{j=0}^5||\pa_x^ju||_{L^2(Q_\e)}\le C(F+M^{1-\tau}F^{\tau}),
\ee
where $\ds M\=\sum_{j=0}^5||\pa_x^ju||_{L^2(Q)}$ and
\bel{23101002}\ba{ll}
\ns\ds F=||(\pa_x^4+\pa_t^{1/3}\pa_x^2+\pa_t^{1/3}\pa_t^{1/3})f(t,x)||_{L^2(Q)}+||f(0,\cdot)||_{L^2(0,L)} \\
\ns\ds\q\q +||(\pa_x^2+\pa_t^{1/3})f(0,x)||_{L^2(0,L)}+\sum_{j=0}^5||h_j(t)||_{H^1(0,T)}.
\ea\ee
\et
\bc
$u\in C^{1,6}(\cl{Q})$ satisfies (\ref{23092101}) with $f=0$ and $h_1=h_2=0$ in $Q$. Then $u=0$ in $Q$.
\ec

We recall the following known result:
\bl\label{23092103} (\cite[Lemma 2.2]{XCY})
If $z\in AC([0,T])$ and  $z(0)=\pa_t^{\g_1}z(0)=0$,
then
$$\pa_t^{\g_2}\pa_t^{\g_1}z=\pa_t^{\g_1+\g_2}z  \mbox{ in } (0,T), \q 0<\g_1+\g_2\le 1, \g_1, \g_2>0.$$
\el

\ \\
{\bf Proof of Theorem \ref{thmapp}. } We divide the proof into several steps.
\ \\
{\bf Step 1. } Set
\bel{23101101}
\hat{u}(t,x)=u(t,x)-\frac{3}{\G(1/3)}t^{1/3}f(0,x)-\frac{3/2}{\G(2/3)}t^{2/3}f_1(x)
\ee
where
\bel{1212-01}
f_1(x)\=(\pa_x^2+\pa_t^{1/3})f(0,x).
\ee

For $\g_1,\g_2\in(0,1)$, we know that
$$
\pa_t^{\g_1} t^{\g_2}=\frac{\G(1+\g_2)}{\G(1+\g_2-\g_1)}t^{\g_2-\g_1}.
$$
Hence, by (\ref{23092101}) and (\ref{23101101}), it is easy to see that
\bel{23101105}
\hat{u}(0,x)=0,\q \pa_t^{1/3}\hat{u}(0,x)=0.
\ee
Further, by (\ref{23101105}) and (\ref{1212-01}), we have
\bel{23101106}\ba{ll}
\ns\ds\pa_t^{2/3}\hat{u}(0,x)=\pa_t^{1/3}[\pa_t^{1/3}\hat{u}(t,x)]\big|_{t=0} \\
\ns\ds =\pa_t^{1/3}\[\pa_x^2u(t,x)+f(t,x)-f(0,x)-\frac{3}{\G(1/3)}t^{1/3}f_1(x)\]\Big|_{t=0} \\
\ns\ds =\pa_x^2\[\pa_x^2u(t,x)+f(t,x)\]\Big|_{t=0}+\pa_t^{1/3}f(0,x)-f_1(x) \\
\ns\ds =(\pa_x^2+\pa_t^{1/3})f(0,x)-f_1(x)=0.
\ea\ee
Next, by (\ref{23101105})--(\ref{23101106}), we have
\bel{23101108}\ba{ll}
\ns\ds\pa_t\hat{u}(t,x)=\pa_t^{1/3}[\pa_t^{2/3}\hat{u}(t,x)] \\
\ns\ds =\pa_t^{1/3}\[\pa_x^4u(t,x)+\pa_x^2f(t,x)+\pa_t^{1/3}f(t,x)-f_1(x)\] \\
\ns\ds =\pa_x^4[\pa_x^2u(t,x)+f(t,x)]+\pa_t^{1/3}(\pa_x^2+\pa_t^{1/3})f(t,x)\\
\ns\ds =\pa_x^6u(t,x)+(\pa_x^4+\pa_t^{1/3}\pa_x^2+\pa_t^{1/3}\pa_t^{1/3})f(t,x).
\ea\ee

Combining (\ref{23101101}) and (\ref{23101108}), we end up with
\bel{23101109}\ba{ll}\ds
\pa_tu-\pa_x^6u&\ds=\wt{f}(t,x),
\ea\ee
where
\bel{1212-02}\ba{ll}\ds
\wt{f}(t,x)=&\ds\frac{1}{\G(1/3)}t^{-2/3}f(0,x)+\frac{3/2}{\G(2/3)}t^{-1/3}f_1\\
\ns&\ds+(\pa_x^4+\pa_t^{1/3}\pa_x^2+\pa_t^{1/3}\pa_t^{1/3})f.
\ea\ee

By (\ref{23101110}) and (\ref{23101109}), we have
\bel{23101003}\left\{\ba{ll}
\ns\ds \pa_tu-\pa_x^6u=\wt{f}(t,x),&\ (t,x)\in Q, \\
\ns\ds \pa_x^ju(t,x)|_{x=0}=h_j(t),\ j=0,1,...,5,&\ t\in(\delta_0,T).
\ea\right.\ee

\ \\
{\bf Step 2. }
Let $\chi(t,x)\in C^\i(Q;[0,1])$ satisfy
\bel{23101005}
\chi(t,x)=\left\{\ba{ll}
\ns\ds 1,&\ (t,x)\in Q_{\e/2}, \\
\ns\ds 0,&\ (t,x)\in Q\backslash Q_{\e/3}.
\ea\right.\ee
 Then we set
\bel{23101006}
\wt{u}=\chi\[u-\sum_{j=0}^5 \frac{x^j}{j!}h_j(t)\]\=\chi u-\chi\wt{h}(t,x).
\ee

The boundary of $Q_{\e/4}$ consists of two line segments:  one is a subset of $\{(t,x)|\delta_0<t<T,\ x=0\}$, while another one  is in $Q\backslash Q_{\e/3}$. One can easily check that $\pa_t\wt{u}(t,x)=\pa_x^{j}\wt{u}(t,x)=0,\ j=0,1,...,5$ on $\pa Q_{\e/4}$, hence $\wt{u}\in H_0^{1,6}(Q_{\e/4})$. Moreover
\bel{23101007}\ba{ll}
\ns\ds\pa_t\wt{u}+\pa_x^6\wt{u} \\
\ns\ds =\chi\wt{f}+(\pa_t\chi)u+\sum_{j=0}^5C_6^j(\pa_x^{6-j}\chi)\pa_x^ju-(\pa_t\chi+\pa_x^6\chi)\wt{h}-\chi\pa_t\wt{h}.\\
\ea\ee

Applying the Carleman estimate in Theorem \ref{thmSW00} with $v=\wt{u}$, $w=e^{\l\psi}\wt{u}$ , we conclude that there is  a large $\l>\l_0>0$ on $Q_{\e/4}$, such that
\bel{23101008}\ba{ll}
&\ds \sum_{j=0}^5||e^{\l\psi}\pa_x^j\wt{u}||^2_{L^2(Q_{\e/4})} \le C||e^{\l\psi}(\pa_t+\pa_x^6)\wt{u}||^2_{L^2(Q_{\e/4})} \\
\ns&\ds \le C||e^{\l\psi}\chi\wt{f}||^2_{L^2(Q_{\e/4})}+C||e^{\l\psi}(\pa_t\chi)u||^2_{L^2(Q_{\e/4})}\\
\ns&\ds\q +C\sum_{j=0}^5||e^{\l\psi}(\pa_x^{6-j}\chi)\pa_x^ju||^2_{L^2(Q_{\e/4})}+C\sum_{j=0}^5||h_j||^2_{H^1(\delta_0,T)}.
\ea\ee

By the definition of $\chi(t,x)$ in (\ref{23101006}), one knows that $\pa_t\chi,\pa_x\chi,...,\pa_x^5\chi$ are supported on $\{x|\psi\le \e/2\}$.  Recall that $\ds M=\sum_{j=0}^5||\pa_x^ju||_{L^2(Q)}$. Thus, we have
\bel{23101009}
||e^{\l\psi}(\pa_t\chi)u||^2_{L^2(Q_{\e/4})}+\sum_{j=0}^5||e^{\l\psi}(\pa_x^{6-j}\chi)\pa_x^ju||^2_{L^2(Q_{\e/4})}\le Ce^{\e\l}M^2.
\ee

Therefore
\bel{23101010}
\sum_{j=0}^5||e^{\l\psi}\pa_x^j\wt{u}||^2_{L^2(Q_{\e/4})}\le Ce^{C\l}(||\wt{f}||^2_{L^2(Q)}+\sum_{j=0}^5||h_j||^2_{H^1(\delta_0,T)})+Ce^{\e\l}M^2.
\ee

Since $u=\wt{u}+\wt{h}$, (\ref{23101010}) yields
\bel{23101011}\ba{ll}
&\ds\sum_{j=0}^5||e^{\l\psi}\pa_x^ju||^2_{L^2(Q_{\e})}\\
\ns&\ds\le C\sum_{j=0}^5||e^{\l\psi}\pa_x^j\wt{u}||^2_{L^2(Q_{\e})}+C||e^{\l\psi}\wt{h}||^2_{L^2(Q_\e)} \\
\ns&\ds\le C\sum_{j=0}^5||e^{\l\psi}\pa_x^j\wt{u}||^2_{L^2(Q_{\e/4})}+Ce^{2\e\l}\sum_{j=0}^5||h_j||^2_{L^2(\delta_0,T)} \\
\ns&\ds\le Ce^{C\l}\[||\wt{f}||^2_{L^2(Q)}+\sum_{j=0}^5||h_j||^2_{H^1(\delta_0,T)}\]+Ce^{\e\l}M^2.
\ea\ee

Meanwhile, it follows from (\ref{23101006}) and (\ref{23101002}) that
\be
||\wt{f}||^2_{L^2(Q)}+\sum_{j=0}^5||h_j||^2_{H^1(\delta_0,T)}\le CF^2.
\eee
Since $\psi>\e$ in $Q_{\e}$, the above, along with (\ref{23101011}), implies
\bel{23101012}
\sum_{j=0}^5||\pa_x^ju||^2_{L^2(Q_{\e})}\le e^{-2\e\l}\sum_{j=0}^5||e^{\l\psi}\pa_x^ju||^2_{L^2(Q_{\e})}\le Ce^{C\l}F^2+Ce^{-\e\l}M^2,
\ee
for all $\l>\l_1$, where $\l_1$ is a sufficiently large number. Replacing the constant $C$ by $Ce^{C\l_1}$, setting $\wt{\l}\=\l-\l_1>0$,  we have
\bel{23101013}\ba{ll}\ds
\sum_{j=0}^5||\pa_x^ju||^2_{L^2(Q_{\e})}&\ds\le Ce^{C\l_1}\[Ce^{C\wt{\l}}F^2+Ce^{-\e\l-C\l_1}M^2\]\\
\ns&\ds\le Ce^{C\wt{\l}}F^2+Ce^{-\e\wt{\l}}M^2.
\ea\ee

\ \\
{\bf Step 3. }
In the case that  $F<M$, we put $\ds
\wt{\l}=\frac{2}{\e+C}\ln{\frac{M}{F}}>0.
$
Then we have
\bel{23101015}
\sum_{j=0}^5||\pa_x^ju||^2_{L^2(Q_{\e})}\le CM^{\frac{2C}{\e+C}}F^{\frac{2\e}{\e+C}},\mbox{ with }\tau\={\frac{2\e}{\e+C}}.
\ee
In the case that $F\ge M$, it follows from (\ref{23101013}) that
\bel{23101016}
\sum_{j=0}^5||\pa_x^ju||^2_{L^2(Q_{\e})}\le C(e^{C\wt{\l}}+e^{-\e\wt{\l}})F^2\le CF^2.
\ee
Combining (\ref{23101015}) and (\ref{23101016}), we get the desired result. \endpf

\end{document}